\documentclass{lmcs}


\usepackage{amsthm}
\usepackage{amssymb}
\usepackage{braket}
\usepackage{prettyref}
\usepackage[hidelinks]{hyperref}

\newcommand{\paren}[1]{\left(#1\right)}
\newcommand{\abs}[1]{\left|#1\right|}
\newcommand{\sqbra}[1]{\left[#1\right]}

\newcommand{\eqw}{\mathrel{=_{w}}}
\newcommand{\neqw}{\mathrel{\neq_{w}}}
\newcommand{\tow}{\mathrel{\rightarrow_{w}}}
\newcommand{\toli}{\mathrel{\rightarrow_{\mathrm{li}}}}
\newcommand{\tor}{\mathrel{\Rightarrow}}
\newcommand{\tailsim}{\mathrel{\sim_{\mathrm{tail}}}}

\newcommand{\Invariant}[1]{\mathrm{Tr}_{#1}}
\newcommand{\ShiftTr}{\mathcal{S}}

\newcommand{\Nat}{\mathbb{N}}
\newcommand{\Cons}{\mathcal{C}}
\newcommand{\Vars}{\mathcal{V}}
\newcommand{\Terms}{\mathcal{T}}
\newcommand{\Alg}{\mathcal{A}}
\newcommand{\NF}{\mathcal{N}}

\newrefformat{cor}{Corollary~\ref{#1}}

\begin{document}

\title[Bluebirds and mockingbirds]{Bluebirds and mockingbirds cannot produce a fixed-point combinator}

\author[T.~Imamura]{Takuma Imamura\lmcsorcid{0000-0001-6491-9137}}

\address{Independent Researcher, Kyoto, Japan}
\email{imamura.takuma.66s@kyoto-u.jp}

\begin{abstract}
Let $B$ be the bluebird combinator with reduction rule $Bxyz \to_{w} x\left(yz\right)$, let $M$ be the mockingbird combinator with reduction rule $Mx \to_{w} xx$, and let $I$ be the identity bird combinator with reduction rule $Ix \to_{w} x$.
A fixed-point combinator, called a sage bird by Smullyan, is a closed term $Y$ such that, for a fresh variable $x$, $Yx$ is equivalent to $x\left(Yx\right)$ under these reduction rules.
For a fixed variable $x$, we construct an invariant $\mathrm{Tr}_{x}\left(u\right)$ of a $BMI$--term $u$ with respect to $\to_{w}$.
This invariant traces the occurrences of $x$ in the leftmost-innermost reduction sequence of $u$.
We then prove that $\mathrm{Tr}_{x}\left(Yx\right) \neq \mathrm{Tr}_{x}\left(x^{r}\left( Yx \right)\right)$ for every $x$--free $BMI$--term $Y$ and every $r\geq 1$.
Consequently, there exists no fixed-point combinator in $BMI$--combinatory logic.
This provides a negative answer to the problem posed by Smullyan in 1985.
\end{abstract}

\maketitle

\section*{Introduction}
\label{sec:introduction}

Smullyan~\cite{Smu00} asked whether a fixed-point combinator $Y$ can be constructed from the bluebird combinator $B$ and the mockingbird combinator $M$ under the following reduction rules:
\begin{align*}
Bxyz &\tow x\paren{yz}, \\
Mx &\tow xx.
\end{align*}
Here, $Y$ is called a fixed-point combinator if $Yx \eqw x\paren{Yx}$ holds for a fresh variable $x$, where $\eqw$ denotes the equivalence relation induced by weak reduction $\tow$.
This problem is also known as Problem~4 in the \emph{TLCA List of Open Problems}~\cite{TLCA} and has remained unsolved since 1985.

This combinatory logic lies on the boundary between the existence and non-existence of fixed-point combinators.
For example, one can construct a fixed point of each individual term.
For every $x$, put
\[
F_{x} := M\paren{BxM}.
\]
Then
\[
F_{x} = M\paren{BxM} \tow BxM\paren{BxM} \tow x\paren{M\paren{BxM}} = xF_{x}.
\]
Thus, the obstruction is not the existence of fixed points, but their \emph{uniform} selection by a combinator.
If the $\lambda$--abstraction $\lambda x. F_{x}$ could be represented in $BM$--combinatory logic, this would yield a fixed-point combinator.
In bases containing a permutative combinator such as $Cxyz \tow xzy$, such $\lambda$--abstractions are possible:
\begin{align*}
\lambda x. M\paren{BxM} &\rightsquigarrow BM \paren{\lambda x. BxM} \\
&\rightsquigarrow BM \paren{CBM}.
\end{align*}
The resulting term actually gives a fixed-point combinator:
\[
BM\paren{CBM} x \tow M \paren{CBMx} \tow M \paren{BxM} = F_{x}.
\]
However, no such permutative term exists in $BM$--combinatory logic.
This is the `permutative effect' that Statman~\cite{Sta11} pointed out when discussing the problem.

There is, moreover, a closed $BM$--term which comes remarkably close to being a fixed-point combinator.
Let
\begin{align*}
T := BM\paren{B\paren{BM}B}.
\end{align*}
The term $Tx$ has the same B\"ohm tree as the application $Yx$ of a fixed-point combinator $Y$ to $x$~\cite{Sta93, Sta11}:
for every $n\geq 1$, there is a residual term $R_{n}$ (with variable $x$) such that
\[
Tx\tow^{\ast} x^{n} R_{n},
\]
where $x^{n} t = x\paren{x\paren{\cdots x\paren{t}\cdots}}$.
However, since the residual term $R_{n}$ is not weakly equivalent to $Tx$, the term $T$ is not a fixed-point combinator.

The \emph{reduction} version of Smullyan's problem was settled earlier.
This problem asks whether there exists a $Y$ such that $Yx \tow^{\ast} x\paren{Yx}$.
Indeed, McCune and Wos~\cite{MW91} and Statman~\cite{Sta93} proved the absence of such a $Y$.
However, this result does not answer the original problem concerning the existence of a $Y$ with $Yx \eqw x(Yx)$.

In an unpublished manuscript~\cite{Sta11}, Statman proposed a negative solution to the original problem.
However, this proof contained a flaw.
Suppose to the contrary that there exists a fixed-point combinator $Y$ consisting of $B$ and $M$.
Statman first showed that $Y$ must be $\beta\eta$--equivalent to a term in a specific class called `polynomials'\footnote{Terms of the form $B^{n} B$, $B^{n} M$ or $B^{n} I$ are called monomials. Terms obtained from monomials by repeating the operation $t \circ u := Btu$ are called polynomials. For example, the term $T := BM\paren{B\paren{BM}B}$ is a polynomial, because $T = M\circ \paren{\paren{BM} \circ B}$.}. (This part is valid.)
He then showed that, in a fragment of Girard's system $F_{1}$, all polynomials have a type of the form $\tau_{1} \to \tau_{2}$, and therefore, $Yx$ has a $\beta\eta$--normal form.
However, there is a counterexample to the latter part.
Indeed, while the aforementioned $T$ is a polynomial, $Tx$ has no $\beta\eta$--normal form.
Therefore, this problem remains open \cite{Sta17}.

The purpose of the present paper is to solve Smullyan's problem.
In fact, we prove the following stronger statement: if $Y$ is a $BMI$--term and $x$ is a fresh variable, then, for every $r\geq 1$,
\[
Y x \neqw x^{r} \paren{Yx}.
\]
In particular, taking $r=1$, no closed $BMI$--term is a fixed-point combinator.
This provides a negative answer to Smullyan's problem.

To distinguish $Yx$ and $x^{r}\paren{Yx}$, we first introduce an invariant $\Invariant{x}\paren{u}$ of a $BMI$--term $u$ with respect to $\tow$.
This invariant traces the occurrences of $x$ in the reduction sequence of $u$ with respect to leftmost-innermost reduction $\toli$.
More precisely, at each $M$--contraction
\[
\Gamma\sqbra{MA} \tow \Gamma\sqbra{AA},
\]
the invariant $\Invariant{x}\paren{u}$ keeps track of two integers: the number of occurrences of $x$ lying strictly before the redex $MA$, and the number of occurrences of $x$ in the duplicated argument $A$.
Since an infinite $\toli$--reduction sequence $u \toli u_{1} \toli u_{2} \toli \cdots$ must contain an infinite number of $M$--contractions, it determines an infinite sequence of pairs of integers:
\[
\paren{p_{0}, q_{0}}, \paren{p_{1}, q_{1}}, \ldots.
\]
The invariant $\Invariant{x}\paren{u}$ is then defined as the equivalence class of the sequence $\paren{p_{i}, q_{i}}$ under a suitable equivalence relation.
On the other hand, if there is no such infinite $\toli$--reduction sequence, $u$ has a weak normal form $N$, and $\Invariant{x}\paren{u}$ is defined to be $N$.
Although $\Invariant{x}$ is defined using leftmost-innermost reduction, it is invariant under arbitrary weak reductions; see \prettyref{thm:invariance}.

We then prove that for every $x$--free $BMI$--term $Y$ and every $r\geq 1$,
\[
\Invariant{x}\paren{Yx} \neq \Invariant{x}\paren{x^{r} \paren{Yx}}.
\]
Since $\Invariant{x}$ is invariant under weak reduction, we conclude that
\[
Y x \neqw x^{r} \paren{Yx}.
\]
As an immediate corollary, we prove the absence of fixed-point combinators in $BMI$--combinatory logic.

This paper is organised as follows.
In \prettyref{sec:preliminaries}, we recall terminology and definitions of combinatory logic.
In \prettyref{sec:invariant}, we introduce the invariant $\Invariant{x}$ and show that it is invariant under weak reduction.
In \prettyref{sec:main}, we prove the main results.We show that $\Invariant{x}$ separates $Yx$ and $x^{r}\paren{Yx}$ for every $x$--free $BMI$--term $Y$ and every $r\geq 1$.
In \prettyref{sec:examples}, we compute the invariant $\Invariant{x}$ for several terms, thereby illustrating the idea that forms the basis of our proof.
Finally, in \prettyref{sec:conclusion}, we conclude by clarifying the scope of our results and stating the open problem concerning fixed-point combinators under $\beta\eta$--equivalence.

\section{Preliminaries}
\label{sec:preliminaries}

\subsection{Terms}

Let $\Cons := \set{B,C,I,M,W}$ be the set of constants used in this paper.
Let $\Vars$ be a set of variables disjoint from $\Cons$.
For a basis $\mathcal{B} \subseteq \Cons$, the set $\Terms_{\mathcal{B}}$ of terms is defined by the following BNF:
\[
u ::= c \mid x \mid \paren{uu} \qquad \paren{c\in \mathcal{B},\ x\in\Vars}.
\]
We abbreviate a basis by omitting braces and commas; for example, $BM$, $BMI$, and $BWI$ denote the bases $\set{B,M}$, $\set{B,M,I}$, and $\set{B,W,I}$, respectively.

Throughout the paper, we regard the application $\paren{u, v} \mapsto \paren{uv}$ as left-associative ($uvw = \paren{uv}w$) and omit parentheses wherever possible.

A term is said to be \emph{$x$--free} if it contains no occurrence of a variable $x$, and is said to be \emph{closed} if it contains no variables.

For $x \in \Vars$ and $r \in \Nat := \set{0, 1, 2, \ldots}$, the term $x^{r} u$ is defined by induction as follows:
\[
x^{0} u := u, \qquad x^{r+1} u := x\paren{x^{r} u}.
\]

For a variable $x$ and a term $u$, let $\abs{u}_x$ denote the number of occurrences of $x$ in $u$.

In what follows, we identify each term with its binary application tree (syntax tree).
The leaves are ordered from left to right.

\subsection{Weak reduction and leftmost-innermost reduction}

We first define the \emph{contraction relation} $\tor$ on $\Terms_{\Cons}$ as follows:
\begin{align*}
Buvw &\tor u\paren{vw}, & Cuvw &\tor uwv, \\
Ix &\tor x, & Mu &\tor uu, \\
Wxy &\tor xyy.
\end{align*}
Given a basis $\mathcal{B} \subseteq \Cons$, the \emph{weak reduction relation} $\tow$ is the least relation on $\Terms_{\mathcal{B}}$ generated by the following rules:
\begin{align*}
u \in\Terms_{\mathcal{B}},\ v\in\Terms_{\mathcal{B}},\ u \tor v &\implies u \tow v, \\
u \tow u' &\implies uv \tow u'v, \\
v \tow v' &\implies uv \tow uv'.
\end{align*}
The reflexive transitive closure of $\tow$ is denoted by $\tow^{\ast}$ and the equivalence closure of $\tow$ by $\eqw$.

In a one-step reduction $\Gamma\sqbra{t} \tow \Gamma\sqbra{u}$ arising from a contraction $t \tor u$, $t$ is called the \emph{weak redex} and $u$ its \emph{contractum}.
A term is said to be \emph{weakly normal} if it contains no weak redex.
A \emph{weak normal form} of a term $u$ is a weak normal term $N$ such that $u \tow^{\ast} N$.
For the sake of brevity, we omit `weak(ly)' if there is no risk of confusion.
Let
\[
\NF_{\mathcal{B}} := \set{u \in \Terms_{\mathcal{B}} | u\text{ is normal}}.
\]

The \emph{leftmost-innermost reduction relation} $\toli$ is then defined as the least relation satisfying the following rules:
\begin{align*}
u \in \NF_{\mathcal{B}},\ v \in \NF_{\mathcal{B}},\ uv \tor w &\implies uv \toli w, \\
u \toli u' &\implies uv \toli u'v, \\
u \in \NF_{\mathcal{B}},\ v\toli v' &\implies uv \toli uv'.
\end{align*}
Its reflexive transitive closure is denoted by $\toli^{\ast}$.

If a term $u$ is not normal, there exists exactly one term $v$ such that $u \toli v$.
Otherwise, there is no term $v$ such that $u \toli v$.
Consequently, every term $u$ generates a unique maximal $\toli$--sequence starting from $u$.

\begin{lem}\label{lem:li-normalisation}
For all $t \in \Terms_{BMI}$ and $N \in \NF_{BMI}$,
\[
t \tow^{\ast} N \iff  t \toli^{\ast} N.
\]
In particular, the $\toli$--sequence of a $BMI$--term terminates if and only if the term has a weak normal form.
\end{lem}

\begin{proof}
(See also Klop~\cite[Section 2.1]{Klo92}.)
Since $\toli$ implies $\tow$, the right-to-left implication is immediate.
To prove the forward direction, we first observe that $BMI$--combinatory logic is a non-erasing, orthogonal term rewriting system.

Now, suppose that $t \tow^{\ast} N$.
By Church's uniform normalisation theorem~\cite{KOO01}, $t$ is strongly normalisable, that is, every $\tow$--sequence from $t$ terminates.
In particular, the maximal $\toli$--sequence from $t$ terminates at some $P \in \NF_{BMI}$.
Moreover, orthogonal term rewriting systems are confluent~\cite{Ros73}, so $t \tow^{\ast} N$ and $t \tow^{\ast} P$ imply $P=N$.
Therefore $t \toli^{\ast} N$.
\end{proof}

For readers familiar with the $\lambda$--calculus, the standard $\lambda$--interpretations of $BMI$--terms lie in the $\lambda_{I}$--calculus, and every contraction is simulated by $\beta$--reduction.
The $\lambda_{I}$--calculus permits an abstraction $\lambda x.u$ only when at least one $x$ occurs free in $u$, so $\beta$--contraction never erases its argument.
Every $\lambda_{I}$--term having a $\beta$--normal form is strongly $\beta$--normalising, a property analogous to that used in \prettyref{lem:li-normalisation}.

\section{An invariant}\label{sec:invariant}

In this section, we introduce an invariant $\Invariant{x}\paren{u}$ of a term $u \in \Terms_{BMI}$ and prove its fundamental properties.

\subsection{Traces}

We first define the $x$--trace $\tau_{x}\paren{u}$ for a term $u$ with no normal form, and the equivalence relation $\tailsim$ on $x$--traces.

\begin{lem}\label{lem:infinite-events}
Let $u_{0} \toli u_{1} \toli \cdots$ be an infinite $\toli$--sequence in $\Terms_{BMI}$.
The $M$--contraction $\Gamma\sqbra{Mv} \toli \Gamma\sqbra{vv}$ appears in the sequence infinitely many times.
\end{lem}

\begin{proof}
Otherwise, there would exist an infinite reduction sequence consisting solely of contractions other than $M$.
However, since contractions other than $M$ strictly reduce the total number of occurrences of constants and variables, this would contradict the assumption.
\end{proof}

\begin{defi}[$x$--trace]
Let $u_{0} \toli u_{1} \toli \cdots$ be an infinite $\toli$--sequence in $\Terms_{BMI}$ and let $x \in \Vars$.
For $i \in \Nat$, let
\[
u_{m_{i}} = \Gamma_{i}\sqbra{M v_{i}} \toli \Gamma_{i}\sqbra{v_{i}v_{i}} = u_{m_{i}+1}
\]
be the $i$--th contraction of $M$ in the sequence.
Define the \emph{$x$--trace} $\tau_{x}\paren{u_{0}}$ of $u_{0}$ by
\[
\tau_{x}\paren{ u_{0} } := \paren{\paren{p_{i}, q_{i}}}_{i\in \Nat} \in \paren{\Nat^{2}}^{\Nat},
\]
where
\begin{align*}
p_{i} &:= \text{the number of occurrences of $x$ in $u_{m_{i}}$ strictly before the redex $M v_{i}$}, \\
q_{i} &:= \abs{v_{i}}_{x}.
\end{align*}
\end{defi}

\subsection{Trace invariants}

Next, we present the definition of the invariant $\Invariant{x}\paren{u}$ of a term $u$.

\begin{defi}[Tail equivalence]
For $\alpha,\beta \in \paren{\Nat^{2}}^{\Nat}$, define \emph{tail equivalence} by
\[
\alpha \tailsim \beta :\iff \exists i, j \in \Nat.\  \forall k \in \Nat.\ \alpha\paren{i + k} = \beta\paren{j + k}.
\]
\end{defi}

Note that, since $i$ and $j$ can be chosen independently in the definition of tail-equivalence, this differs from the eventual equality of sequences.

\begin{defi}[Trace invariants]
Define the set $\Alg_{BMI}$ as the disjoint union:
\[
\Alg_{BMI} := \NF_{BMI} \sqcup \paren{\Nat^{2}}^{\Nat} / \tailsim.
\]
For a term $u \in \Terms_{BMI}$ and a variable $x$, define the \emph{trace invariant} $\Invariant{x}\paren{u} \in \Alg_{BMI}$ as follows:
\[
\Invariant{x}\paren{u} := \begin{cases}
    N, & \text{$u$ has a normal form $N$}, \\
    \sqbra{\tau_{x}\paren{u}}, & \text{otherwise}.
\end{cases}
\]
\end{defi}

\begin{lem}\label{lem:prefix}
If $u \toli^{\ast} v$, then $\Invariant{x}\paren{u}=\Invariant{x}\paren{v}$.
\end{lem}

\begin{proof}
If $v$ has a normal form $N$, then so does $u$.
Hence $\Invariant{x}\paren{u}=\Invariant{x}\paren{v} = N$.

Suppose that $v$ has no normal form.
Let $u \toli w_{1} \toli \cdots \toli w_{n} = v$ be a $\toli$--sequence from $u$ to $v$.
Let $c$ be the number of $M$--contractions in that sequence.
Then $\tau_{x}\paren{u}_{i+c} = \tau_{x}\paren{v}_{i}$ for all $i\in \Nat$, so $\tau_{x}\paren{u} \tailsim \tau_{x}\paren{v}$.
Therefore $\Invariant{x} \paren{u} = \sqbra{\tau_{x}\paren{u}} = \sqbra{\tau_{x}\paren{v}}=\Invariant{x}\paren{v}$.
\end{proof}

\begin{defi}[Shift operators]
For $c \in \Nat$, define the map
\[
S_{c}\colon \paren{\Nat^{2}}^{\Nat} \to \paren{\Nat^{2}}^{\Nat}
\]
by 
\[
S_{c}\alpha :=\paren{\paren{p_{i} + c, q_{i}}}_{i\in\Nat},
\]
where
\[
\alpha = \paren{\paren{p_{i}, q_{i}}}_{i\in\Nat}.
\]
Since $S_{c}$ preserves tail equivalence, it induces a map
\[
\ShiftTr_{c} \colon \paren{\Nat^{2}}^{\Nat}/\tailsim \to \paren{\Nat^{2}}^{\Nat}/\tailsim,\qquad \sqbra{\alpha} \mapsto \sqbra{S_{c}\alpha}.
\]
\end{defi}

\subsection{Algebraic structure of trace invariants}

We now define a binary operator $\star_{x}$ on $\Alg_{BMI}$, making it an algebraic semantics for $BMI$--combinatory logic.

\begin{defi}
Let $x \in \Vars$. 
Define $\star_{x} \colon \Alg_{BMI} \times \Alg_{BMI} \to \Alg_{BMI}$ by
\[
a \star_{x} b := \begin{cases}
    a, & a \in \paren{\Nat^{2}}^{\Nat} / \tailsim, \\
    \ShiftTr_{\abs{a}_{x}} b & a \in \NF_{BMI},\ b \in \paren{\Nat^{2}}^{\Nat} / \tailsim, \\
    \Invariant{x}\paren{ab}, & a,b \in \NF_{BMI}.
\end{cases}
\]
\end{defi}

\begin{lem}
Let $x \in \Vars$.
The operator $\star_{x} \colon \Alg_{BMI} \times \Alg_{BMI} \to \Alg_{BMI}$ is well-defined.
\end{lem}

\begin{proof}
We consider three cases, corresponding to the three clauses in the definition of $\star_{x}$.

\begin{enumerate}
\item $a \in \paren{\Nat^{2}}^{\Nat} / \tailsim$.
The value $a \star_{x} b = a$ is an element of $\Alg_{BMI}$ and involves no choice of representatives.

\item $a \in \NF_{BMI}$ and $b \in \paren{\Nat^{2}}^{\Nat} / \tailsim$.
Suppose that $b = \sqbra{\alpha} = \sqbra{\beta}$.
Put $c := \abs{a}_{x}$.
Since $\alpha \tailsim \beta$, there exist $i,j \in \Nat$ such that $\alpha\paren{i+k} = \beta\paren{j+k}$ for every $k \in \Nat$.
Adding $c$ to the first coordinate gives
\[
\paren{S_{c}\alpha}\paren{i+k} = \paren{S_{c}\beta}\paren{j+k} \qquad \paren{k \in \Nat}.
\]
Thus $S_{c}\alpha \tailsim S_{c}\beta$, so $a \star_{x} b = \ShiftTr_{c}b$ is independent of the representative of $b$.

\item $a, b \in \NF_{BMI}$.
This case involves no choice of representatives. \qedhere
\end{enumerate}
\end{proof}

\begin{thm}[Homomorphism]\label{thm:homomorphism}
For all $u,v \in \Terms_{BMI}$ and $x\in \Vars$,
\[
\Invariant{x}\paren{uv} = \Invariant{x}\paren{u} \star_{x} \Invariant{x}\paren{v}.
\]
\end{thm}

\begin{proof}
We consider three cases according to whether $u$ and $v$ have normal forms.

\begin{enumerate}
\item The $\toli$--sequence of $u$ is infinite.
Write this sequence as $u = u_{0} \toli u_{1} \toli \cdots$.
By the definition of $\toli$, the $\toli$--sequence of $uv$ is
\[
uv = u_{0}v \toli u_{1}v \toli \cdots.
\]
Each $M$--contraction $\Gamma\sqbra{Mw}\toli \Gamma\sqbra{ww}$ in $u$'s sequence corresponds precisely to the $M$--contraction $\Gamma\sqbra{Mw}v\toli \Gamma\sqbra{ww}v$ in $uv$'s sequence, and vice versa.
So
\[
\Invariant{x}\paren{uv} = \sqbra{\tau_{x}\paren{u}} = \Invariant{x}\paren{u} = \Invariant{x}\paren{u} \star_{x} \Invariant{x}\paren{v}.
\]

\item $u \toli^{\ast} U$ for some $U \in \NF_{BMI}$, and the $\toli$--sequence of $v$ is infinite.
Write the latter sequence as $v = v_{0} \toli v_{1} \toli \cdots$.
Since $U$ is normal, and each $v_{i}$ is non-normal, the leftmost-innermost strategy selects a redex inside $v_{i}$, rather than at the root of $Uv_{i}$.
Thus, the $\toli$--sequence of $uv$ has the form
\[
uv \toli^{\ast} Uv = Uv_{0} \toli Uv_{1} \toli \cdots.
\]
Each $M$--contraction $\Gamma\sqbra{Mw} \toli \Gamma\sqbra{ww}$ in $v$'s sequence corresponds precisely to the $M$--contraction
\[
U\paren{\Gamma\sqbra{Mw}} \toli U\paren{\Gamma\sqbra{ww}}
\]
in $Uv$'s sequence, and vice versa.
The duplicated argument is $w$ in both contractions, while the latter has exactly $\abs{U}_{x}$ additional occurrences of $x$ before the redex.
Thus a trace pair $\paren{p,q}$ in the former becomes $\paren{p+\abs{U}_{x},q}$ in the latter, so
\[
\tau_{x}\paren{Uv} = S_{\abs{U}_{x}}\tau_{x}\paren{v}.
\]
Since $\Invariant{x}\paren{u} = U$, \prettyref{lem:prefix} gives
\begin{align*}
\Invariant{x}\paren{uv}
&= \Invariant{x}\paren{Uv} \\
&= \sqbra{S_{\abs{U}_{x}}\tau_{x}\paren{v}} \\
&= \ShiftTr_{\abs{U}_{x}}\paren{\Invariant{x}\paren{v}} \\
&= \Invariant{x}\paren{u} \star_{x} \Invariant{x}\paren{v}.
\end{align*}

\item $u \toli^{\ast} U$ and $v \toli^{\ast} V$ for some $U,V \in \NF_{BMI}$.
Write these finite sequences as
\[
u = u_{0} \toli \cdots \toli u_{m} = U,
\qquad
v = v_{0} \toli \cdots \toli v_{n} = V.
\]
Each contraction $\Gamma\sqbra{R} \toli \Gamma\sqbra{R'}$ in the first sequence occurs in $uv$'s sequence as
\[
\Gamma\sqbra{R}v \toli \Gamma\sqbra{R'}v.
\]
Once $U$ is reached, its normality ensures that each contraction $\Delta\sqbra{R} \toli \Delta\sqbra{R'}$ in the second sequence occurs as
\[
U\paren{\Delta\sqbra{R}} \toli U\paren{\Delta\sqbra{R'}}.
\]
Consequently, the $\toli$--sequence of $uv$ begins with
\begin{align*}
uv = u_{0}v
&\toli \cdots \toli u_{m}v = Uv_{0} \\
&\toli \cdots \toli Uv_{n} = UV.
\end{align*}
Therefore, by \prettyref{lem:prefix} and the definition of $\star_{x}$,
\[
\Invariant{x}\paren{uv}
= \Invariant{x}\paren{UV}
= U \star_{x} V
= \Invariant{x}\paren{u} \star_{x} \Invariant{x}\paren{v}.
\]
This applies whether or not $UV$ has a normal form. \qedhere
\end{enumerate}
\end{proof}

\begin{thm}[Constants]\label{thm:identities}
Let $x \in \Vars$.
Regard the constant terms $B,M,I$ as elements of $\Alg_{BMI}$.
Then, for all $a,b,c \in \Alg_{BMI}$,
\begin{itemize}
\item $\paren{\paren{B\star_{x} a} \star_{x} b} \star_{x} c = a \star_{x} \paren{b \star_{x} c}$;
\item $M \star_{x} a = a \star_{x} a$;
\item $I \star_{x} a = a$.
\end{itemize}
\end{thm}

\begin{proof}
To prove the $B$--identity, we consider four cases.
\begin{enumerate}

\item $a \in \paren{\Nat^{2}}^{\Nat}/\tailsim$.
Since $\abs{B}_{x}=0$, we have $B\star_{x}a=\ShiftTr_{0}a=a$.
The first clause in the definition of $\star_{x}$ therefore gives
\[
\paren{\paren{B\star_{x}a}\star_{x}b}\star_{x}c
= a
= a\star_{x}\paren{b\star_{x}c}.
\]

\item $a \in \NF_{BMI}$ and $b \in \paren{\Nat^{2}}^{\Nat}/\tailsim$.
The term $Ba$ is normal, so $B\star_{x}a=Ba$ and $\abs{Ba}_{x}=\abs{a}_{x}$.
Hence
\begin{align*}
\paren{\paren{B\star_{x}a}\star_{x}b}\star_{x}c
&= \paren{\ShiftTr_{\abs{a}_{x}}b}\star_{x}c \\
&= \ShiftTr_{\abs{a}_{x}}b \\
&= a\star_{x}\paren{b\star_{x}c},
\end{align*}
where the last two equalities use the fact that a tail-equivalence class is unchanged by application on the right.

\item $a,b \in \NF_{BMI}$ and $c \in \paren{\Nat^{2}}^{\Nat}/\tailsim$.
The terms $Ba$ and $Bab$ are normal, and
$\abs{Bab}_{x}=\abs{a}_{x}+\abs{b}_{x}$.
Consequently,
\begin{align*}
\paren{\paren{B\star_{x}a}\star_{x}b}\star_{x}c
&= \ShiftTr_{\abs{a}_{x}+\abs{b}_{x}}c \\
&= \ShiftTr_{\abs{a}_{x}}\paren{\ShiftTr_{\abs{b}_{x}}c} \\
&= a\star_{x}\paren{b\star_{x}c}.
\end{align*}
Here we used $\ShiftTr_{r+s}=\ShiftTr_{r}\circ\ShiftTr_{s}$, which follows directly from the definition of the shift operators.

\item $a,b,c \in \NF_{BMI}$.
Since all proper subterms of $Babc$ are normal, its selected contraction is
\[
Babc \toli a\paren{bc}.
\]
By \prettyref{lem:prefix} and \prettyref{thm:homomorphism},
\begin{align*}
\paren{\paren{B\star_{x}a}\star_{x}b}\star_{x}c
&= \Invariant{x}\paren{Babc} \\
&= \Invariant{x}\paren{a\paren{bc}} \\
&= a\star_{x}\paren{b\star_{x}c}.
\end{align*}
\end{enumerate}

For the $M$-- and $I$--identities, first suppose that
$a \in \paren{\Nat^{2}}^{\Nat}/\tailsim$.
Since $\abs{M}_{x}=\abs{I}_{x}=0$, the definition gives
\[
M\star_{x}a=\ShiftTr_{0}a=a=a\star_{x}a,
\qquad
I\star_{x}a=\ShiftTr_{0}a=a.
\]
If $a \in \NF_{BMI}$, the selected contractions are
\[
Ma \toli aa,
\qquad
Ia \toli a.
\]
Thus \prettyref{lem:prefix} and the definition of $\star_{x}$ give
\begin{align*}
M\star_{x}a
&= \Invariant{x}\paren{Ma}
= \Invariant{x}\paren{aa}
= a\star_{x}a, \\
I\star_{x}a
&= \Invariant{x}\paren{Ia}
= \Invariant{x}\paren{a}
= a. \qedhere
\end{align*}
\end{proof}

\begin{thm}[Invariance under weak reduction]\label{thm:invariance}
Let $x \in \Vars$. 
For all $u,v \in \Terms_{BMI}$,
\[
u \eqw v \implies \Invariant{x}\paren{u}=\Invariant{x}\paren{v}.
\]
\end{thm}
\begin{proof}
Immediate from \prettyref{thm:homomorphism} and \prettyref{thm:identities}.
\end{proof}

It should be noted that while the trace invariant $\Invariant{x}$ is defined using a particular leftmost-innermost reduction sequence, it is invariant under an arbitrary weak reduction.

\section{Main result}\label{sec:main}

Having made the above preparations, we now prove the main result of this paper:

\begin{thm}\label{thm:main}
Let $x \in \Vars$ and $r \geq 1$.
For every finite $x$--free term $Y \in \Terms_{BMI}$,
\[
Yx \neqw x^{r}\paren{Yx}.
\]
In particular, $BMI$--combinatory logic has no fixed-point combinator under weak equivalence $\eqw$.
\end{thm}

\subsection{Some structural lemmas}
\label{sec:structural-lemmas}

Fix $x \in \Vars$ and an $x$--free term $P \in \NF_{BMI}$.
The following lemmas concern the $\toli$--sequence starting from $Px$.

\begin{lem}\label{lem:rightmost}
For every $u$ such that $Px \toli^{\ast} u$, the rightmost leaf of $u$ is labelled $x$, and every redex occurrence in $u$ contains that leaf.
\end{lem}

\begin{proof}
We argue by induction on the number of $\toli$--steps.
The assertion holds for $Px$: its rightmost leaf is $x$, and every subterm not containing this leaf is a subterm of the normal term $P$.

Suppose that the assertion holds for $u$ and that
\[
u = \Gamma\sqbra{R} \toli \Gamma\sqbra{R'}=u'.
\]
Since $R$ is a redex and therefore, not normal, it contains the rightmost leaf of $u$ by the induction hypothesis.
Moreover, every proper subterm of $R$ is normal because the contraction is innermost.
There are three cases.
\begin{enumerate}

\item The contraction is
\[
\Gamma\sqbra{Babc} \toli \Gamma\sqbra{a\paren{bc}}.
\]
The rightmost leaf lies in $c$, which is retained unchanged.
The only new subterms appearing on the contractum are $bc$ and $a(bc)$, both of which contain this leaf.
Every other subterm within the contractum is a subterm of one of the normal terms $a,b,c$.

\item The contraction is
\[
\Gamma\sqbra{Ma} \toli \Gamma\sqbra{a_{L}a_{R}},
\]
where $a_{L}$ and $a_{R}$ denote the two occurrences of the same normal term $a$.
The rightmost leaf of $a$ is labelled $x$, so the new rightmost leaf is the final $x$ in $a_{R}$.
Both copies of $a$ are normal, and the new subterm $a_{L}a_{R}$ contains the new rightmost leaf.

\item The contraction is
\[
\Gamma\sqbra{Ia} \toli \Gamma\sqbra{a}.
\]
The normal term $a$ is retained unchanged and contains the rightmost leaf.
\end{enumerate}
In each case, subterms disjoint from the contracted occurrence are unchanged, while every ancestor of that occurrence still contains the rightmost leaf.
Thus the assertion also holds for $u'$.
\end{proof}

Consequently, every redex in a term of this $\toli$--sequence lies on the path from the root to the rightmost leaf.
This path is obtained by repeatedly taking the right child of an application.

\begin{lem}\label{lem:increase}
If $MA, MA'$ are the redexes at two consecutive $M$--contractions in the $\toli$--sequence of $Px$, then
\[
\abs{A'}_{x} \geq \abs{A}_{x}+1.
\]
\end{lem}

\begin{proof}
Write the first of the two $M$--contractions as
\[
\Gamma\sqbra{MA} \toli \Gamma\sqbra{AA}.
\]
Since the contraction is innermost, $A$ is a normal form.
By \prettyref{lem:rightmost}, the redex $MA$ contains the rightmost leaf of the whole term, which is labelled $x$.
Thus the rightmost leaf of $A$ is labelled $x$.

Only $B$-- and $I$--contractions occur between the two consecutive $M$--contractions.
We prove, by induction on the number of these intervening contractions, that every term between the consecutive $M$--contractions has the form
\[
U_{1}\paren{U_{2}\paren{\cdots\paren{U_{k}A}\cdots}},
\]
where $k\geq 1$, all the $U_{i}$ are normal, and the rightmost leaf of $U_{k}$ is labelled $x$.

For the initial term $\Gamma\sqbra{AA}$, observe that the path from the root of $\Gamma\sqbra{MA}$ to the redex $MA$ takes only right-child edges, since $MA$ contains the rightmost leaf of $\Gamma\sqbra{MA}$.
Each left subtree along this path is normal, because $MA$ is the leftmost-innermost redex of $\Gamma\sqbra{MA}$.
Taking these left subterms in order from the root to $MA$ as $U_{1},\ldots,U_{m}$, we have
\[
\Gamma\sqbra{-} = U_{1}\paren{U_{2}\paren{\cdots\paren{U_{m} {-}}\cdots}},
\]
where $m=0$ means $\Gamma\sqbra{-} = {-}$.
Therefore,
\[
\Gamma\sqbra{AA} = U_{1}\paren{U_{2}\paren{\cdots\paren{U_{m}\paren{AA}}\cdots}}.
\]
Setting $k := m+1$ and $U_{k} := A$ gives the required form.

For the induction step, suppose that the current term has the form:
\[
U_{1}\paren{U_{2}\paren{\cdots\paren{U_{k}A}\cdots}}.
\]
Every subterm of the current term is either a subterm of one of $U_{1},\ldots,U_{k},A$ or has the form
\[
U_{i}\paren{\cdots\paren{U_{k}A}\cdots} \qquad \paren{1\leq i\leq k}.
\]
Since $U_{1},\ldots,U_{k},A$ are all normal, the redex that is contracted at this step must have the latter form.
Write this redex as $U_{j}T$.

For a $B$--contraction, we have $U_{j} = Bab$, and the contraction replaces $BabT$ with $a\paren{bT}$.
Thus the factor $U_{j}$ is replaced with the two factors $a,b$.
The new factors $a$ and $b$ are normal, since $U_{j}$ is normal.
If $j < k$, the last factor $U_{k}$ is unchanged.
If $j = k$, the new last factor is $b$, whose rightmost leaf
is the same as the rightmost leaf of $Bab=U_{k}$ and is therefore labelled $x$.
Hence the required form is preserved.

For an $I$--contraction, we have $U_{j} = I$, and the contraction replaces $IT$ with $T$, just deleting the factor $U_{j}$.
Since the rightmost leaf of $U_{k}$ is labelled $x$, we have
$U_{k} \neq I$ and hence $j < k$.
The last factor is therefore unchanged, and at least one factor remains.
This completes the induction.

At the second $M$--contraction, the selected redex has the form
\[
MA' = U_{j}\paren{\cdots\paren{U_{k}A}\cdots}
\]
for some $1\leq j\leq k$.
Thus $U_{j} = M$.
Since the rightmost leaf of $U_{k}$ is labelled $x$, we have $U_{k}\neq M$ and hence $j < k$.
Consequently,
\[
A' = U_{j+1}\paren{\cdots\paren{U_{k}A}\cdots},
\]
and therefore,
\begin{align*}
\abs{A'}_{x}
&= \abs{A}_{x} + \sum_{i=j+1}^{k}\abs{U_{i}}_{x} \\
&\geq \abs{A}_{x} + \abs{U_k}_{x} \\
&\geq \abs{A}_{x} + 1. \qedhere
\end{align*}
\end{proof}

\begin{cor}\label{cor:increasing}
If $Px$ has an infinite $\toli$--sequence and
\[
\tau_{x}\paren{Px} = \paren{\paren{p_{i}, q_{i}}}_{i \in \Nat},
\]
then $q_{0} < q_{1} < q_{2} < \cdots$.
\end{cor}

\begin{proof}
By \prettyref{lem:infinite-events}, there are infinitely many $M$--contractions.
Let $A_i$ be the argument duplicated at the $i$--th $M$--contraction.
By definition, $q_{i} = \abs{A_{i}}_{x}$.
Applying \prettyref{lem:increase} to each consecutive pair gives
\[
q_{i+1} = \abs{A_{i+1}}_{x} \geq \abs{A_{i}}_{x} + 1 = q_{i} + 1 \qquad \paren{i\in\Nat}.\qedhere
\]
\end{proof}

\begin{lem}\label{lem:separation}
Let $\alpha := \paren{\paren{p_{i}, q_{i}}}_{i\in\Nat} \in \paren{\Nat^{2}}^{\Nat}$ and suppose the second coordinates of $\alpha$ are strictly increasing.
Then
\[
\alpha \not\tailsim S_{r}\alpha
\]
for every $r\in\Nat$ with $r\geq1$.
\end{lem}

\begin{proof}
Suppose that $\alpha\tailsim S_{r}\alpha$.
By the definition of tail equivalence, there exist $i,j\in\Nat$ such that
\[
\paren{p_{i+k}, q_{i+k}} = \paren{p_{j+k} + r, q_{j+k}} \qquad \paren{k\in\Nat}.
\]
Taking $k=0$ in the second coordinate gives $q_{i} = q_{j}$.
Since the $q$--coordinate is strictly increasing, $i = j$, so the first coordinate gives $p_{i} = p_{i} + r$.
This contradicts $r\geq 1$.
\end{proof}

\subsection{Proof of the theorem}
We first prove the following lemma:

\begin{lem}\label{lem:head-xs-shift-trace}
For every $u\in\Terms_{BMI}$ and $r\in\Nat$,
\[
\Invariant{x}\paren{x^{r}u} = \begin{cases}
x^{r}N, & \Invariant{x}\paren{u}=N\in\NF_{BMI},\\
\sqbra{S_r\alpha} = \ShiftTr_{r}\sqbra{\alpha}, & \Invariant{x}\paren{u}=\sqbra{\alpha} \in \paren{\Nat^{2}}^{\Nat}/\tailsim.
\end{cases}
\]
\end{lem}

\begin{proof}
Note that the $\toli$--sequence of $x^{r}u$ is obtained from that of $u$ by placing each term in the context $x^{r}\paren{-}$.

Suppose first that $u\toli^{\ast}N$ with $N$ normal.
Then $x^{r}u\toli^{\ast}x^{r}N$, and $x^{r}N$ is normal.

Now suppose that the $\toli$--sequence of $u$ is infinite, and write it as $u=u_{0}\toli u_{1}\toli\cdots$.
The $\toli$--sequence starting from $x^{r}u$ is
\[
x^{r}u = x^{r}u_{0} \toli x^{r}u_{1} \toli \cdots.
\]
Each contraction $\Gamma\sqbra{MA}\toli\Gamma\sqbra{AA}$ corresponds to
\[
x^{r}\paren{\Gamma\sqbra{MA}} \toli x^{r}\paren{\Gamma\sqbra{AA}}.
\]
The duplicated argument is unchanged and exactly $r$ additional occurrences of $x$ precede the redex.
Thus
\[
\tau_{x}\paren{x^{r}u} = S_{r}\tau_{x}\paren{u},
\]
which proves the second clause.
\end{proof}

\begin{proof}[Proof of \prettyref{thm:main}]
Fix an $x$--free term $Y\in\Terms_{BMI}$ and $r\geq 1$.
By \prettyref{thm:invariance}, it suffices to show that
\[
\Invariant{x}\paren{Yx} \neq \Invariant{x}\paren{x^{r}\paren{Yx}}.
\]
Every reduct of an $x$--free term is $x$--free, since the $BMI$--contraction rules introduce no new variables.
There are three cases:
\begin{enumerate}

\item The $\toli$--sequence of $Y$ is infinite.
At each $M$--contraction in this sequence, both the duplicated argument and the whole term are $x$--free.
Hence
\[
\tau_{x}\paren{Y} = \zeta_{0},
\]
where
\[
\zeta_{s} := \paren{\paren{s,0}}_{i\in\Nat}.
\]
By \prettyref{thm:homomorphism} and the first clause in the definition of $\star_x$,
\[
\Invariant{x}\paren{Yx} =\sqbra{\zeta_{0}}\star_{x} x =\sqbra{\zeta_{0}}.
\]
On the other hand, \prettyref{lem:head-xs-shift-trace} gives
\[
\Invariant{x}\paren{x^{r}\paren{Yx}} = \sqbra{S_{r}\zeta_{0}} =\sqbra{\zeta_{r}}.
\]
Since $r\geq 1$, no entry of $\zeta_{0}$ equals an entry of $\zeta_{r}$, so $\zeta_{0} \not\tailsim \zeta_{r}$.
Therefore, $\Invariant{x}\paren{Yx} \neq \Invariant{x}\paren{x^{r}\paren{Yx}}$.

\item $Y\toli^{\ast}P$ and $Px\toli^{\ast}N$ for some $P,N\in\NF_{BMI}$.
The term $P$ is $x$--free, and
\[
Yx \toli^{\ast} Px \toli^{\ast} N.
\]
By \prettyref{lem:prefix} and \prettyref{lem:head-xs-shift-trace},
\[
\Invariant{x}\paren{Yx}=N, \qquad
\Invariant{x}\paren{x^{r}\paren{Yx}}=x^{r}N.
\]
These are distinct normal-form terms, because
\[
\abs{x^{r}N}_{x} = \abs{N}_{x} + r > \abs{N}_{x}.
\]

\item $Y\toli^{\ast}P$ for some $P\in\NF_{BMI}$, and the $\toli$--sequence of $Px$ is infinite.
The term $P$ is $x$--free.
Put $\alpha := \tau_{x}\paren{Px}$.
Since $Yx \toli^{\ast} Px$, \prettyref{lem:prefix} and \prettyref{lem:head-xs-shift-trace} give
\[
\Invariant{x}\paren{Yx}=\sqbra{\alpha}, \qquad
\Invariant{x}\paren{x^{r}\paren{Yx}}=\sqbra{S_{r}\alpha}.
\]
The second coordinates of $\alpha$ are strictly increasing by \prettyref{cor:increasing}.
Therefore $\alpha\not\tailsim S_{r}\alpha$ by \prettyref{lem:separation}, so the two invariants are distinct.
\end{enumerate}
Thus the required inequality holds in every case.
\end{proof}

\begin{cor}\label{cor:BM}
Let $x \in \Vars$ and $r \geq 1$.
For every finite $x$--free term $Y \in \Terms_{BM}$,
\[
Yx \neqw x^{r}\paren{Yx}.
\]
In particular, $BM$--combinatory logic has no fixed-point combinator under weak equivalence $\eqw$.
\end{cor}

\section{Examples}
\label{sec:examples}

In this section, we compute $\Invariant{x}$ for several concrete terms to illustrate the ideas underlying the proof of the main theorem.

\subsection{Non-uniform fixed points}
\label{ex:nonuniform}

The absence of a fixed-point combinator does not exclude fixed points whose construction depends on the given term.
For a variable $x$, put
\[
A_{x} := BxM, \qquad
F_{x} := MA_{x} = M\paren{BxM}.
\]
The term $A_{x}$ is normal and $\abs{A_{x}}_{x} = 1$.
The $\toli$--sequence of $F_{x}$ begins with
\[
F_{x} = MA_{x} \toli A_{x}A_{x} = BxMA_{x} \toli x\paren{MA_{x}} = xF_{x}.
\]
The same $\toli$--reduction occurs inside every outer context $x^{n}\paren{-}$:
\[
x^{n}F_{x} \toli x^{n}\paren{A_{x}A_{x}} \toli x^{n+1}F_{x} \qquad \paren{n\in\Nat}.
\]
This is an infinite $\toli$--sequence.
In its $n$--th $M$--contraction, the duplicated argument is always $A_{x}$, and precisely $n$ occurrences of $x$ precede the redex.
Hence
\[
\tau_{x}\paren{F_{x}} = \beta,
\]
where
\[
\beta := \paren{\paren{i,1}}_{i\in\Nat}.
\]
For every $r \in \Nat$,
\[
S_{r}\beta = \paren{\paren{i+r, 1}}_{i\in\Nat} = \paren{\beta_{i+r}}_{i\in \Nat},
\]
so $\beta\tailsim S_{r}\beta$.

Consequently, \prettyref{lem:head-xs-shift-trace} gives
\[
\Invariant{x}\paren{x^{r}F_{x}} =\sqbra{S_{r}\beta} =\sqbra{\beta} =\Invariant{x}\paren{F_{x}},
\]
which is consistent with $F_{x}\toli^{\ast}x^{r}F_{x}$ and \prettyref{lem:prefix}.

Since $F_{x}$ is not of the form $Px$ with $x$--free normal term $P$, the results in \prettyref{sec:structural-lemmas} cannot apply to this case.
Indeed, the rightmost leaf of $F_{x}$ is $M$, and the second coordinate of $\tau_{x}\paren{F_{x}}$ is constant.

\subsection{A term similar to a fixed-point combinator}
\label{ex:nonstandard}

Consider the following closed normal term (Statman~\cite{Sta93,Sta11}):
\[
t := BM\paren{B\paren{BM}B}.
\]
As we shall see, $tx$ reduces to terms of the form $x^{n}R_{n}$
for arbitrarily large $n$.
In that sense, $t$ behaves similarly to a fixed-point combinator; nevertheless, $tx \neqw x\paren{tx}$.

Let $x$ be a fixed variable.
For $n \in \Nat$, define terms $A_{n}$ by
\begin{align*}
A_{0} &:= BM\paren{Bx}, \\
A_{n+1} &:= \paren{Bx}A_{n},
\end{align*}
i.e., $A_{n} = \paren{Bx}^{n} \paren{BM\paren{Bx}}$.
Clearly, each $A_{n}$ is normal and has exactly $n + 1$ occurrences of the variable $x$.

The $\toli$--sequence of $tx$ starts with
\begin{align*}
tx &= BM\paren{B\paren{BM}B} x \\ 
&\toli M\paren{B\paren{BM}Bx} \\
&\toli M\paren{BM\paren{Bx}} = MA_{0} \\
&\toli A_{0} A_{0} = BM\paren{Bx}A_{0} \\
&\toli M\paren{\paren{Bx}A_{0}} = MA_{1} = x^{0}\paren{MA_{1}}.
\end{align*}
For each $n \geq 1$,
\begin{align*}
MA_{n} &\toli A_{n}A_{n} = \paren{Bx}^{n}A_{0} A_{n} \\
&\toli^{n} x^{n} \paren{A_{0}A_{n}} = x^{n} \paren{BM\paren{Bx}A_{n}} \\
&\toli x^{n} \paren{M\paren{BxA_{n}}} = x^{n}\paren{MA_{n+1}},
\end{align*}
where the first $\toli$ is an $M$--contraction and the other $n+1$ are $B$--contractions.
The $n$--th $M$--contraction of the $\toli$--sequence is therefore
\[
x^{d_{n}}\paren{MA_{n}} \toli x^{d_{n}}\paren{A_{n} A_{n}},
\]
where
\[
d_{n} := \sum_{k=0}^{n-1} k = \frac{n\paren{n-1}}{2}.
\]
It follows that
\[
\tau_{x}\paren{tx} = \paren{\paren{d_{n}, \abs{A_{n}}_{x}}}_{n\in\Nat} = \paren{\paren{d_{n}, n+1}}_{n\in\Nat}.
\]
On the other hand, for each $r \geq 1$,
\[
\tau_{x}\paren{x^{r}\paren{tx}} =  \paren{\paren{d_{n} + r, n+1}}_{n\in\Nat}
\]
by \prettyref{lem:head-xs-shift-trace}.
These traces are not tail-equivalent by \prettyref{lem:separation}, so
\[
\Invariant{x}\paren{tx} \neq \Invariant{x}\paren{x^{r}\paren{tx}}\qquad\paren{r\geq 1}.
\]
By \prettyref{thm:invariance}, we conclude that
\[
tx \neqw x^{r}\paren{tx}\qquad\paren{r\geq 1}.
\]

This calculation also illustrates \prettyref{lem:rightmost} and \prettyref{lem:increase}.
After the $M$--contractions
\[
MA_{n} \toli A_{n}A_{n},
\]
the following $B$--contractions retain the right-hand copy of $A_{n}$ inside the next duplicated argument $A_{n+1} = BxA_{n}$.
The additional occurrence of $x$ gives
\[
\abs{A_{n+1}}_{x} = \abs{A_{n}}_{x} + 1.
\]

\subsection{A fixed-point combinator in the BWI-system}
\label{ex:bwi}

In this subsection, we explain why the analogue of \prettyref{thm:main} does not hold in the $BWI$--system.
The proof of non-existence of fixed-point combinators relies on the structural lemmas in \prettyref{sec:structural-lemmas}, which fail in $BWI$--combinatory logic.

Consider the closed normal term \cite[p.~178]{CFC58}:
\[
Y := B\paren{WI}\paren{BWB}.
\]
Put $A_{x} := W\paren{Bx}$.
Then $A_{x}$ is normal and $\abs{A_{x}}_{x} = 1$.
The $\toli$--sequence of $Yx$ begins with
\begin{align*}
Yx
&\toli \paren{WI}\paren{BWBx} \\
&\toli \paren{WI}\paren{W\paren{Bx}} = WIA_{x}\\
&\toli \paren{IA_{x}}A_{x} \toli A_{x}A_{x}.
\end{align*}
Since $A_{x}$ is normal, reduction continues as
\begin{equation}\label{eq:bwi-cycle}\tag{$\dagger$}
A_{x}A_{x}=W\paren{Bx}A_{x} \toli BxA_{x}A_{x} \toli x\paren{A_{x}A_{x}}.
\end{equation}
The two steps contract $W$ and $B$, respectively.
Consequently,
\[
Yx\eqw A_{x}A_{x}\eqw x\paren{A_{x}A_{x}}\eqw x\paren{Yx},
\]
so $Y$ is a fixed-point combinator in the $BWI$--system.

For comparison, define the $x$--trace $\tau_{x}^{BWI}$ for the $BWI$--system.
At each $W$--contraction
\[
\Gamma\sqbra{Wuv} \toli \Gamma\sqbra{\paren{uv}v},
\]
consider the pair $\paren{p,\abs{v}_{x}}$, where $p$ counts the occurrences of $x$ preceding the redex $Wuv$.
The analogues of \prettyref{lem:li-normalisation} and \prettyref{lem:infinite-events} hold for the $BWI$--system by the same arguments, with $W$--contractions replacing $M$--contractions in the latter.
We may therefore define, just as in the $BMI$--case,
\[
\Invariant{x}^{BWI}\paren{u} := \begin{cases}
N, & u\toli^{\ast}N\in\NF_{BWI},\\
\sqbra{\tau_{x}^{BWI}\paren{u}}, & \text{otherwise}.
\end{cases}
\]

The corresponding algebra is
\[
\Alg_{BWI} := \NF_{BWI}\sqcup\paren{\Nat^{2}}^{\Nat}/\tailsim,
\]
with the operator $\star_{x}^{BWI}$ defined in the same way.
The same arguments establish the analogues of
\prettyref{thm:homomorphism}, \prettyref{thm:identities}, and \prettyref{thm:invariance}, with the $M$--identity replaced by
\[
\paren{W\star_{x}^{BWI}a}\star_{x}^{BWI}b = \paren{a\star_x^{BWI}b}\star_{x}^{BWI}b.
\]

The $\toli$--reduction \eqref{eq:bwi-cycle} occurs inside every outer context $x^{n}\paren{\cdot}$.
Its $W$--contraction duplicates $A_{x}$ and has exactly $n$ occurrences of $x$ before the redex.
Hence
\[
\tau_x^{BWI}\paren{A_{x}A_{x}} =\beta,
\]
where
\[
\beta := \paren{\paren{i,1}}_{i\in\Nat}.
\]
The initial $W$--contraction $WIA_{x} \toli \paren{IA_{x}}A_{x}$ contributes one additional pair $\paren{0, 1}$.
There are no other $W$--contractions before $A_{x}A_{x}$ is reached, so
\[
\tau_{x}^{BWI}\paren{Yx}
= \paren{\paren{0,1}, \paren{0,1}, \paren{1,1}, \paren{2,1}, \ldots}.
\]
Deleting its first entry gives $\beta$, and $S_1\beta$ is obtained from $\beta$ by deleting its first entry.
Therefore
\[
\Invariant{x}^{BWI}\paren{Yx} =\sqbra{\beta} =\sqbra{S_{1}\beta} =\Invariant{x}^{BWI}\paren{x\paren{Yx}}.
\]
The last equality follows from the same outer-context trace calculation used in \prettyref{lem:head-xs-shift-trace}.

The failure of the two structural lemmas can be seen directly in these contractions.
First, although $u,v$ are normal in an innermost $W$--redex $Wuv$, the new left subterm $uv$ of the $W$--contractum $\paren{uv}v$ need not be normal.
In particular,
\[
WIA_{x}\toli\paren{IA_{x}}A_{x}
\]
creates the redex $IA_{x}$, which does not contain the rightmost leaf of the whole term.
Thus the analogue of \prettyref{lem:rightmost} for the $BWI$--system fails.

Second, the analogue of \prettyref{lem:increase} for $W$--contractions fails.
In the proof of \prettyref{lem:increase}, the argument duplicated by the next $M$--contraction contains both the previously duplicated argument $A$ and a factor $U_k$ whose rightmost leaf is labelled $x$.
It therefore contains at least $\abs{A}_x+1$ occurrences of $x$.
By contrast, in the present example,
\[
W\paren{Bx}A_{x} \toli \paren{BxA_{x}}A_{x} \toli x\paren{A_{x}A_{x}} = x\paren{W\paren{Bx}A_{x}}.
\]
Every subsequent $W$--contraction duplicates the same term $A_{x}$, so the number of occurrences of $x$ in the duplicated argument remains equal to $1$.

\section{Conclusion}\label{sec:conclusion}

We constructed a trace invariant $\Invariant{x}$ of a term in $BMI$--combinatory logic.
Using this invariant, we proved that, for any $r\geq 1$, there exists no $x$--free term $Y$ such that $Yx \eqw x^{r}\paren{Yx}$.
In particular, the absence of a fixed-point combinator was obtained as the case where $r = 1$.

Our results concern weak equivalence and do not settle the corresponding questions for $\beta$-- or $\beta\eta$--equivalence of the $\lambda$--interpretations.
For example, the present invariant is not compatible with $\beta$--equivalence of the $\lambda$--interpretations.
Indeed, under translation to $\lambda$--terms, $BII$ is $\beta$--equivalent to $I$:
\[
BII = \paren{\lambda xyz.x\paren{yz}} II \to_{\beta}^{\ast} \lambda z.I\paren{Iz} \to_{\beta}^{\ast} \lambda z.z = I.
\]
However, their invariants do not coincide:
\[
\Invariant{x}\paren{BII} = BII \neq I = \Invariant{x}\paren{I}.
\]

The following problem, which Statman~\cite{Sta11} attempted to solve, remains unresolved:

\begin{prob}
Is there an applicative combination $Y$ of $B:=\lambda xyz.x\paren{yz}$, $M:=\lambda x.xx$ and $I:=\lambda x.x$ such that $Yx =_{\beta\eta} x\paren{Yx}$?
\end{prob}

It might be worth attempting to distinguish between $Yx$ and $x\paren{Yx}$ by constructing an appropriate $\beta\eta$--invariant, just as was done for weak equivalence in this paper.

\section*{Acknowledgements}

The author is grateful to the reviewers for their valuable comments, which helped improve the clarity of the manuscript.

\bibliographystyle{alphaurl}
\bibliography{ref}

@book{CFC58,
    author = {Haskell H. Curry and Robert Feys and William Craig},
    title = {Combinatory {L}ogic},
    bookseries = {Studies in Logic and the Foundations of Mathematics},
    publisher = {North-Holland Publishing},
    address   = {Amsterdam},
    year = {1958},
    issn = {0049-237X},
    doi = {https://doi.org/10.1016/S0049-237X(08)72031-7},
}

@inproceedings{KOO01,
    author    = {Khasidashvili, Zurab and Ogawa, Mizuhito and Oostrom, Vincent van},
    title     = {Uniform Normalisation beyond Orthogonality},
    year      = {2001},
    publisher = {Springer-Verlag},
    address   = {Berlin, Heidelberg},
    booktitle = {Proceedings of the 12th International Conference on Rewriting Techniques and Applications},
    pages     = {122--136},
    numpages  = {15},
    series    = {RTA '01}
}

@incollection{Klo92,
    author      = {Klop, J. W.},
    title       = {Term rewriting systems},
    year        = {1992},
    isbn        = {0198537611},
    publisher   = {Oxford University Press},
    address     = {New York, NY, USA},
    booktitle   = {Handbook of Logic in Computer Science},
    volume      = {2},
    subtitle    = {Background: Computational Structures},
    pages       = {1--116},
    numpages    = {116},
}

@article{MW91,
    author  = {William McCune and Larry Wos},
    title   = {The absence and the presence of fixed point combinators},
    journal = {Theoretical Computer Science},
    volume  = {87},
    number  = {1},
    pages   = {221--228},
    year    = {1991},
    doi     = {https://doi.org/10.1016/0304-3975(91)90034-Y},
}

@article{Ros73,
    author     = {Rosen, Barry K.},
    title      = {Tree-Manipulating Systems and {C}hurch-{R}osser Theorems},
    year       = {1973},
    issue_date = {Jan. 1973},
    publisher  = {Association for Computing Machinery},
    volume = {20},
    number = {1},
    doi = {10.1145/321738.321750},
    journal = {Journal of the Association for Computing Machinery},
    month = jan,
    pages = {160--187},
    numpages = {28},
}

@book{Smu00,
    author    = {Raymond Smullyan},
    title     = {To {M}ock a {M}ockingbird: and {O}ther {L}ogic {P}uzzles},
    publisher = {Oxford University Press},
    address   = {Oxford and New York},
    year      = 2000,
    note      = {Reprint of Knopf 1985 edition},
}

@article{Sta93,
    author  = {Richard Statman},
    title   = {Some examples of non-existent combinators},
    journal = {Theoretical Computer Science},
    volume  = {121},
    number  = {1},
    pages   = {441--448},
    year    = {1993},
    doi     = {https://doi.org/10.1016/0304-3975(93)90096-C},
}

@misc{Sta11,
    author = {Richard Statman},
    title  = {To Type A Mockingbird},
    year   = {2011},
    month  = dec,
    howpublished = {https://tlca.di.unito.it/PAPER/TypeMock.pdf},
    note   = {Unpublished manuscript},
}

@incollection{Sta17,
    author    = {Richard Statman},
    title     = {Some Tweets About Mockingbirds},
    editor    = {Fitting, Melvin and Rayman, Brian},
    booktitle = {Raymond Smullyan on Self Reference},
    publisher = {Springer International Publishing},
    address   = {Cham},
    pages     = {113--122},
    year      = {2017},
    isbn      = {978-3-319-68732-2},
    doi       = {10.1007/978-3-319-68732-2_6},
}

@misc{TLCA,
    author = {Richard Statman},
    title  = {Problem \# 4, {TLCA} List of Open Problems},
    howpublished = {https://tlca.di.unito.it/opltlca/},
    note   = {Last retrieved: 2026-09-17},
}

\end{document}